\documentclass[letterpaper, 10 pt, conference]{ieeeconf}  

\IEEEoverridecommandlockouts                              

\usepackage{afterpage}
\usepackage{epsfig} 
\usepackage{times} 
\usepackage{amsmath} 
\usepackage{amssymb}  
\usepackage{amsfonts}
\usepackage{mathptmx} 
\usepackage{tabularx}
\usepackage{epstopdf}
\usepackage{subfigure}
\usepackage{float}
\usepackage{epic,color}
\usepackage{mathrsfs}
\usepackage{dsfont,MnSymbol}
\usepackage{algorithm}
\usepackage{algorithmic}
\usepackage{multirow} 
\usepackage[T1]{fontenc}
\newcounter{rmnum}
\newenvironment{romannum}{\begin{list}{{\upshape (\roman{rmnum})}}{\usecounter{rmnum}
\setlength{\leftmargin}{14pt}
\setlength{\rightmargin}{8pt}
\setlength{\itemsep}{2pt}
\setlength{\itemindent}{-1pt}
}}{\end{list}}

\newcounter{anum}

\newlength{\noteWidth}
\long\def\notes#1{\ifinner
             {\tiny #1}
             \else
              \marginpar{\parbox[t]{\noteWidth}{\raggedright\tiny #1}}
               \fi}

\def\IEEEQEDclosed{\mbox{\rule[0pt]{1.3ex}{1.3ex}}}

\def\qed{\ifmmode\IEEEQEDclosed\else{\unskip\nobreak\hfil
\penalty50\hskip1em\null\nobreak\hfil\IEEEQEDclosed
\parfillskip=0pt\finalhyphendemerits=0\endgraf}\fi}

\def\qed{\hspace*{\fill}~\IEEEQED\par\endtrivlist\unskip}

\def\Re{\mathbb{R}}

\def\notes#1{\marginpar{\tiny #1}\typeout{Notes!
Notes!
Notes!
}}
\renewcommand{\notes}[1]{\typeout{notes!}}

\newcommand{\tr}{\mbox{tr}}

\def\Re{\field{R}}
\def\k{{\sf K}}

\def\clZ{{\cal Z}}

\def\Expect{{\sf E}}

\def\Expect{{\sf E}}

\def\P{{\sf P}}

\def\IEEEQEDclosed{\mbox{\rule[0pt]{1.3ex}{1.3ex}}}
\def\qed{\nobreak\hfill\IEEEQEDclosed}

\def\clZ{{\cal Z}}

\newtheorem{lemma}{Lemma}
\newtheorem{remark}{Remark}
\newtheorem{proposition}{Proposition}

\def\beq{\begin{eqnarray}} 
\def\bc{\begin{center}} 
\def\be{\begin{enumerate}}
\def\bi{\begin{itemize}} 
\def\bs{\begin{small}}
\def\bS{\begin{slide}}
\def\ec{\end{center}} 
\def\ee{\end{enumerate}}
\def\ei{\end{itemize}}
\def\es{\end{small}}
\def\eS{\end{slide}}
\def\eeq{\end{eqnarray}}

\newcommand{\newP}[1]{\medskip \noindent{\bf #1:}}

\newcommand{\PP}{{\sf P}}

\newcommand{\ud}{\,\mathrm{d}}

\def\Re{\mathbb{R}}

\def\Expect{{\sf E}}

\def\clZ{{\cal Z}}

\renewcommand{\Re}{\mathbb{R}}

\newcommand{\X}{X}

\newcommand{\NN}{\mathcal{N}}

\newcommand{\mN}{m^{(N)}}
\newcommand{\SigN}{\Sigma^{(N)}}

\newcommand{\Fnorm}[1]{\|#1\|_F}
\newcommand{\Ricc}{\text{Ricc}}
\newcommand{\sRicc}{\sqrt{\text{Ricc}}}

\newcommand{\fSpace}{\mathcal{S}}

\title{\LARGE \bf
Error Analysis of the Stochastic Linear Feedback Particle Filter}

\author{Amirhossein Taghvaei, Prashant G. Mehta 
\thanks{Financial support from the NSF CMMI grants 1334987 and 1462773 is gratefully acknowledged. 
}
\thanks{A.~Taghvaei and P.~G.~Mehta are with the Coordinated
  Science Laboratory and the Department of Mechanical Science and
  Engineering at the University of Illinois at Urbana-Champaign (UIUC).}
}

\begin{document}
\maketitle
\thispagestyle{empty}
\pagestyle{empty}

\begin{abstract}
This paper is concerned with the convergence and long-term stability analysis of the feedback particle filter (FPF) algorithm. The FPF is an interacting system of $N$ particles where the interaction is designed such that the empirical distribution of the particles approximates the  posterior distribution. It is known that in the mean-field limit ($N=\infty$), the distribution of the particles is equal to the posterior distribution. 
However little is known about the convergence to the mean-field limit. In this paper, we consider the FPF algorithm for the linear Gaussian setting.
In this setting, the algorithm is similar to the ensemble Kalman-Bucy filter algorithm.  
Although these algorithms have been numerically evaluated and widely used in applications, their convergence and long-term stability analysis remains an  active area of research. 
In this paper, we show that, (i) the mean-field limit is well-defined with a unique strong solution; (ii) the mean-field process is stable with respect to the initial condition; (iii) we provide conditions such that the finite-$N$ system is long term stable and we obtain some mean-squared error estimates that are uniform in time.      	
\end{abstract}
\section{Introduction}
Feedback particle filter (FPF) is a numerical algorithm to approximate the solution of the nonlinear filtering problem~\cite{taoyang_TAC12,yang2016}. 
The algorithm is comprised of a system of $N$ interacting particles. The interaction is designed such that the empirical distribution of the particles approximates the posterior distribution. 
The FPF algorithm is an alternative to the sequential importance sampling and resampling particle filters~\cite{doucet09}. The salient feature of the FPF,  compared to the conventional particle filters, is that it replaces the importance sampling and resampling step with a feedback control law. 
Because of this difference, in numerical evaluations, FPF does not suffer from issues such as particle degeneracy that is commonly observed in the conventional particle filters~\cite{doucet09}. 
Also in various numerical evaluations and comparisons, it has been observed that FPF exhibit smaller simulation variance and better scaling properties with the problem dimension compared to particle filters ~\cite{berntorp2015,stano2014,surace_SIAM_Review}.

In the mean-field ($N=\infty$) limit,  the FPF is known to be exact, i.e, the conditional probability distribution of the particles is equal to the posterior distribution. However, little is known about the convergence of the finite-$N$ system to the mean-field limit and its long-term stability. The objective of this paper is to address some of these questions in the linear Gaussian setting. 

In the linear Gaussian setting, the FPF algorithm is similar to the ensemble Kalman filter algorithm~\cite[Sec. 4.3]{yang2016}. 
Ensemble Kalman filter (EnKF) was first introduced in~\cite{evensen1994sequential}, in discrete time setting, as an alternative to the extended Kalman filter (EKF) for applications in geophysical sciences. 
In these applications, the state dimension is typically  very high. The main advantage of the EnKF, compared to the EKF, is that the computational cost of the EnKF scales linearly with the state dimension whereas the computational cost of the EKF scales as the dimension squared. 

Since its introduction, the EnKF has evolved into different formulations. The most two well-known formulations are (i) EnKF based on perturbed observation~\cite{evensen2003ensemble} and (ii) the square root EnKF~\cite{whitaker2002ensemble}. For a review of the different discrete time formulations of the EnKF see~\cite[Ch. 6-7]{jdw:ReichCotter2015}~\cite[Ch. 4]{law2015data}.       

The two aforementioned discrete time formulations of the EnKF algorithm have been extended to the continuous time setting~\cite{jdw:BergemannReich2012}. The continuous time formulation of the EnKF is usually referred to as the ensemble Kalman-Bucy filter (EnKBF). 
For a recent review of the EnKBF algorithm and its connection to the FPF algorithm see~\cite{TaghvaeiASME2017}. The EnKBF algorithm and the linear FPF have the following three established formulations:
\begin{romannum}
	\item EnKBF with perturbed observation~\cite{jdw:BergemannReich2012}~\cite{delmoral2016stability};
	\item Stochastic linear FPF~\cite[Eq. (26)]{yang2016} which is same as  the square root EnKBF~\cite{jdw:BergemannReich2012};
	\item Deterministic linear FPF~\cite[Eq. (15)]{AmirACC2016}~\cite{jana2016stability};  
\end{romannum}

In our previous conference publication~\cite{AmirACC18}, we presented the analysis of the deterministic linear FPF. The objective of this paper is to extend the analysis to the stochastic linear FPF. 
The contributions of the paper are as follows: 
\begin{romannum}
	\item we show that the mean-field limit is well defined, and a unique solution exists~(Prop.~\ref{prop:mean-field-existence});
	\item we show the long-term stability of the mean-field system with respect to initial distribution~(Prop.~\ref{prop:mean-field-stability});
	\item we provide estimates for the mean-square error for any finite $N$ that are uniform in time~(Prop.~\ref{prop:mean-var-conv} and Prop.~\ref{prop:conv-POC}).
\end{romannum}
\newP{Literature review on error analysis of EnKF}
Theoretical error and convergence analysis of the EnKF algorithm is an active area of research.   
In the discrete time setting, it is shown that the ensemble distribution converges to the mean-field limit with the convergence rate $O(\frac{1}{\sqrt{N}})$ for any finite time~\cite{gland2009}~\cite{mandel2011convergence}. The asymptotic (in time) stability analysis is more difficult. It is shown that if the system dynamics is stable and admits a Lyapunov function, and the observation model satisfies the "observable energy criterion" (which holds under full state observation), then the system is ergodic and it is stable with respect to initial conditions~\cite{tong2016nonlinear}. The well-posedness of the EnKF and its accuracy using the variance inflation technique is studied in~\cite{stuart2014stability}. Related  finite-time results on the convergence of the discrete-time square root EnKF appear in~\cite{kwiatkowski2015convergence}. The analysis in~\cite{kwiatkowski2015convergence} is  simpler as the model is deterministic and the update formula exactly equals the Kalman filter update formula.

The analysis for EnKBF and linear FPF is more recent. 
For EnKBF with perturbed observation, under certain assumptions (stable and fully observable), it has been shown that the empirical distribution of the ensemble converges to the mean-field distribution uniformly for all time with the rate $O(\frac{1}{\sqrt{N}})$~\cite{delmoral2016stability}. This result has  been extended to the nonlinear setting for the case with Langevin type dynamics with a strongly convex potential and full linear observation~\cite{delmoral2017stability}.

Analysis of the deterministic linear FPF is easier because the update formuala is identical to the Kalman filter update formula. For the linear Gaussian setting, it is shown that (i) the empirical distribution converges to the mean-field limit for any finite time; (ii) and even for a finite number of particles, the long term error converges to zero~\cite{AmirACC18}. The convergence and long term stability results are shown for the nonlinear setting as well, where it is assumed that drift function is Lipschitz and the system is fully observed with small measurement noise~\cite{jana2016stability}. 

\newP{Notation} For a vector $m$, $|m|$
denotes the Euclidean norm.  
For a square matrix $\Sigma$,
$\Fnorm{\Sigma}$
denotes the Frobenius norm,
$\|\Sigma\|_2$ 
is the spectral norm,
$\Sigma^\top$ is the matrix-transpose, $\tr(\Sigma)$ is the
matrix-trace, and $\text{cond}(\Sigma)= \|\Sigma\|_2 \|\Sigma^{-1}\|_2$ is the condition number.
The space of symmetric positive definite matrices is denoted by $S^d_{++}$. 
$\NN(m,\Sigma)$ denotes a Gaussian probability distribution with mean
$m$ and covariance $\Sigma \in S^d_{++}$.
The $L^2$-Wasserstein distance between  two probability measures $\mu,\nu$ is denoted by $W_2(\mu,\nu)$. For a positive integer $n$, double factorial $n!!=\prod_{k=0}^{\lfloor\frac{n}{2}\rfloor-1}(n-2k)$.

There are three types of stochastic process
considered in this paper: (i) $X_t$ denotes the state of the (hidden)
signal at
time $t$; (ii) ${X}_t^i$ denotes the state of the
$i^{\text{th}}$ particle in a population of $N$ particles; and (iii)
$\bar{X}_t$ denotes the state of the McKean-Vlasov model obtained in
the mean-field limit ($N=\infty$).  
The mean and the covariance for
these are denoted as follows: (i) ($m_t,\Sigma_t$) is the conditional mean and
the conditional covariance pair for $X_t$; (ii) ($\mN_t,\SigN_t$) is the
empirical mean and the empirical covariance for the ensemble
$\{X_t^i\}_{i=1}^N$; and (iii) ($\bar{m}_t,\bar{\Sigma}_t$) is the
conditional mean and the conditional covariance for $\bar{X}_t$.  
The notation is tabulated in
Table~\ref{tab:symbols-states}. 

\begin{table}[h]
	\centering
	\begin{tabular}{|c|c|c|}
		\hline
		Variable & Notation & Equation \\ \hline\hline 
		State of the hidden process & $X_t$ & Eq.~\eqref{eq:dyn}\\		
		State of the mean-field process  & $\bar{X}_t$ & Eq.~\eqref{eq:mean-field-FPF} \\ 
		State of the $i^{\text{th}}$ particle in finite-$N$ sys.& $X_t^i$
		&Eq.~\eqref{eq:finite-N-FPF}
		\\
		i.i.d copies of the mean-field process & $\bar{X}_t^i$
		&Eq.~\eqref{eq:mean-field-copy} 
		\\
		\hline
		\vspace*{-0.1in}\\
		Kalman filter mean and covariance & ${m}_t,{\Sigma}_t$ & Eq.~\eqref{eq:KF-mean}-\eqref{eq:KF-variance} 
		\\
		Mean-field mean and covariance & $\bar{m}_t,\bar{\Sigma}_t$ & Eq.~\eqref{eq:mean-field-mean-cov} 
		\\		
		Empirical mean and covariance & $\mN_t,\SigN_t$ & Eq.~\eqref{eq:empr_app_mean_var}
		\\
		\hline
		\vspace*{-0.1in}\\
		Error process for mean-field system & $\bar{\xi}_t$ & Eq.~\eqref{eq:mean-field-error}
		\\ 
		Error process for finite-$N$ system  & $\xi^i_t$ & Eq.~\eqref{eq:finite-N-error} 
		\\
		Error processes driven by $\sRicc(Q_t)$& $\bar{\xi}^{(Q)}_t$
		&Eq.~\eqref{eq:mean-field-error-Q} 
		\\
		i.i.d error processes for $\bar{X}^i_t$ & $\bar{\xi}_t^i$
		&Eq.~\eqref{eq:mean-field-error-copy} 
		\\
		\hline
		\vspace*{-0.1in}\\
		Process noise for hidden process& $B_t$ & Eq.~\eqref{eq:dyn}
		\\ 
		Process noise for mean-field process  & $\bar{B}_t$ & Eq.~\eqref{eq:mean-field-FPF}
		\\
		Process noise for finite-$N$ system  & $B^i_t$ & Eq.~\eqref{eq:finite-N-FPF}
		\\
		Process noise for the  mean  & $B^{(N)}_t$ & Eq.~\eqref{eq:finite-N-mean}
		\\
		Process noise for the covariance  & $M_t$ & Eq.~\eqref{eq:finite-N-DRE}
		\\
		\hline	
	\end{tabular}
	\caption{Nomenclature}  
	\label{tab:symbols-states}
\end{table}
\section{Problem formulation}
\newP{Linear Gaussian filtering problem} Consider the linear Gaussian filtering problem:
\begin{subequations}
	\begin{align}
	\ud X_t &=  AX_t \ud t + \sigma_B\ud B_t\label{eq:dyn}\\
	\ud Z_t &= HX_t\ud t + \ud W_t\label{eq:obs}
	\end{align}
\end{subequations}
where $X_t \in \Re^d$ is the (hidden) state at time $t$, $Z_t \in \Re^m$ is the
observation; $A$, $H$, $\sigma_B$ are matrices of appropriate
dimension; and $\{B_t\}$, $\{W_t\}$ are mutually independent Wiener
processes taking values in $\Re^{d_B}$ and $\Re^m$, respectively. Without
loss of generality, the covariance matrices associated with $\{B_t\}$
and $\{W_t\}$ are identity matrices.
The initial condition $X_0$ is drawn from a Gaussian
distribution $\mathcal{N}(m_0,\Sigma_0)$, independent of $\{B_t\}$ and
$\{W_t\}$. The filtering problem is to compute the posterior
distribution $\P(X_t|\clZ_t)$ where $\clZ_t:=\sigma(Z_s;s\in[0,t])$
denotes the time-history of observations up to time $t$ (filtration).

\newP{Kalman-Bucy filter} For the linear Gaussian problem~\eqref{eq:dyn}-\eqref{eq:obs}, the posterior distribution
$\PP(\X_t|\clZ_t)$ is Gaussian $\mathcal{N}(m_t,\Sigma_t)$, whose mean
and covariance are given by the Kalman-Bucy filter~\cite{kalman-bucy}:
\begin{subequations}
	\begin{align}
	\ud m_t &= A m_t\ud t + \mathsf{K}_t(\ud Z_t - Cm_t\ud t)\label{eq:KF-mean}\\
	\frac{\ud \Sigma_t}{\ud t} &= 
	A \Sigma_t +  \Sigma_tA^\top + \sigma_B\sigma_B^{\top} - \Sigma_tC^\top C\Sigma_t
	\label{eq:KF-variance}
	\end{align}
\end{subequations}
where $\k_t:=\Sigma_tC^\top$ is the Kalman gain, and the filter is
initialized with the prior $\mathcal{N}(m_0,\Sigma_0)$. 

\newP{FPF algorithm} The main steps of the FPF algorithm are to: (i) construct a stochastic process, denoted by $\bar{X}_t$, whose posterior distribution (given $\clZ_t$) is equal to the posterior distribution of $X_t$; (ii) and then simulate $N$ stochastic process, denoted by $\{X^i_t\}_{i=1}^N$, to empirically approximate the distribution of $\bar{X}_t$. 
\begin{equation*}
\Expect[f(X_t)|\clZ_t]\overset{\text{step  (i)}}{=}\Expect[f(\bar{X}_t)|\clZ_t]\overset{\text{step (ii)}}{\approx} \frac{1}{N}\sum_{i=1}^Nf(X^i_t)
\end{equation*}
The process $\bar{X}_t$ is referred to as mean-field process and the $N$ processes $\{X^i_t\}_{i=1}^N$ are referred to as particles.

\newP{Mean-field process} 
The evolution of $\bar{X}_t$ is given by the sde:
\begin{equation}
\ud \bar{X}_t = A \bar{X}_t \ud t + \sigma_B  \ud \bar{B}_t + \bar{\k}_t(\ud Z_t - \frac{H\bar{X}_t + H\bar{m}_t}{2}\ud t)
\label{eq:mean-field-FPF}
\end{equation} 
where $\bar{B}_t$ is an independent copy of the process noise $B_t$, $\bar{\k}_t:=\bar{\Sigma}_tH^\top$ is the Kalman gain, the mean-field terms 
\begin{equation} \label{eq:mean-field-mean-cov}
\bar{m}_t := \Expect[\bar{X}_t|\clZ_t],\quad \bar{\Sigma}_t := \text{var}(X_t|\clZ_t)
\end{equation}
and the initial condition  $\bar{X}_0 \sim \mathcal{N}(m_0,\Sigma_0)$.

\newP{Finite-$N$ system} The evolution of the particles $\{X_t^i\}_{i=1}^N$ is given by
the sde:
\begin{align}
\ud X^i_t &= A X^i_t\ud t + \sigma_B \ud B_t^i + \k^{(N)}_t (\ud Z_t -
\frac{HX^i_t + Hm^{(N)}_t}{2}\ud t) 
\label{eq:finite-N-FPF}
\end{align}
for $i=1,\ldots,N$ where $\{B^i_t\}_{i=1}^N$ are
independent copies of the process noise $B_t$; $\k^{(N)}=\Sigma^{(N)}_tH^\top$ is the Kalman gain, $X^i_0
\stackrel {\text{i.i.d}}{\sim} \mathcal{N}(m_0,\Sigma_0)$; and 
\begin{align}
m^{(N)}_t&:=\frac{1}{N}\sum_{i=1}^N X^i_t, \quad \Sigma^{(N)}_t
:=\frac{1}{N-1}\sum_{i=1}^N (X^i_t-m^{(N)}_t)(X^i_t-m^{(N)}_t)^\top
\label{eq:empr_app_mean_var}
\end{align}



The sde~\eqref{eq:mean-field-FPF} represents the mean-field limit of the interacting particle system~\eqref{eq:finite-N-FPF}. These models are referred to as McKean-Vlasov SDEs~\cite{mckean1966} and their analysis is referred to as {\it propagation of chaos}~\cite{sznitman1991}.

\newP{Paper outline} 
The objective of this paper is to present the propagation of chaos analysis for the linear FPF model~\eqref{eq:mean-field-FPF}-\eqref{eq:finite-N-FPF}. After presenting necessary background about stability of the Kalman-Bucy filter in Sec.~\ref{sec:stability-KF}, we present analysis of the mean-field model in~Sec.~\ref{sec:mean-field-process} and the convergence of the finite-$N$ system in~Sec.~\ref{sec:finite-N}.

Throughout the paper, we make the following assumption:

\newP{Assumption A1} The pair $(A,H)$ is detectable and $(A,\sigma_B)$ is stabilizable.

\newP{Assumption A2} The covariance matrix $\Sigma_B:=\sigma_B\sigma_B^\top \succ 0$.
\section{Stability of the Kalman filter}\label{sec:stability-KF}
\subsection{Ricatti flow}
%
For the linear Gaussian filtering problem~\eqref{eq:dyn}-\eqref{eq:obs}
 %
 %
 \begin{align*}
 \Ricc(Q) &:= A Q +  Q A^\top + \Sigma_B -
 Q H^\top HQ
 \end{align*}
 for $Q \in S^d_{++}$. Define the algebraic Riccati equation (ARE) and the differential Riccati equation (DRE) as follows:
 \begin{subequations}
 	\begin{align}
 \text{(DRE)}:& \quad 	\frac{\ud \Sigma_t}{\ud t} = \Ricc(\Sigma_t)\label{eq:DRE}
 	\\
 	\text{(ARE)}:& \quad \Ricc(\Sigma)=0\label{eq:ARE}
 	\end{align}	
 \end{subequations}
Let $\Phi_{t,s}$  be the state transition matrix for the linear time-varying flow
\begin{equation}
\frac{\ud}{\ud t}\Phi_{t,s} = (A-\Sigma_t H^\top H)\Phi_{t,s},\quad \Phi_{s,s}=I\label{eq:Riccati-flow}
\end{equation}
where $\Sigma_t$ is the solution to the DRE~\eqref{eq:DRE} with initial condition $\Sigma_0$ at $t=0$.
\medskip


 \begin{lemma} \label{lem:Riccati} Consider the ARE~\eqref{eq:ARE}, DRE~\eqref{eq:DRE}, and the flow~\eqref{eq:Riccati-flow}. Then 
 	\begin{romannum} 
 		\item ~\cite[Sec. 23 Thm. 1]{brockett2015finite} There exists a unique positive definite solution $\Sigma_\infty$ to the ARE~\eqref{eq:ARE}.
 		\item The explicit solution to the DRE~\eqref{eq:DRE} is given by:
 		\begin{equation}
 		\Sigma_t = \Sigma_\infty + e^{F_\infty t}D_t^{-1}e^{F_\infty^\top t}\label{eq:DRE-sol}
 		\end{equation}
 		where $F_\infty:=A - \Sigma_\infty H^\top H$, and  $D_t:=(\Sigma_0-\Sigma_\infty)^{-1} + \int_0^t e^{F_\infty^\top s}H^\top H e^{F_\infty s}\ud s$.  
 		\item\cite[Sec. 23, Thm. 3]{brockett2015finite} The eigenvalues of $F_\infty$ have strictly negative real part, i.e,  
 		\begin{equation*}
 		\lambda_0:=\min \{-\text{real}(\lambda):~\lambda~\text{is an eigenvalue of $F_\infty$}\}>0
 		\end{equation*} 
%
 		\item \cite[Eq. (16)]{ocone1996} For all $\lambda<\lambda_0$, there exists constant $\kappa>0$ and a time $t_0>0$ such that 
 		 \begin{equation*}
 		 \|\Phi_{t,s}\|_2\leq \kappa e^{-\lambda(t-s)},\quad \forall t\geq s \geq t_0
 		 \end{equation*}

 	\end{romannum}
 \end{lemma}
\medskip
 		
%
%

	It follows from Lemma~\ref{lem:Riccati} that the Kalman filter is stable in the following sense: Let $(m_t,\Sigma_t)$ and $(\tilde{m}_t,\tilde{\Sigma}_t)$ be solutions to the Kalman filter equations~\eqref{eq:KF-mean}-\eqref{eq:KF-variance} starting from different initial conditions $(m_0,\Sigma_0)$ and $(\tilde{m}_0,\tilde{\Sigma}_0)$ respectively. Then for all $\lambda<\lambda_0$, there are constants $M_1,M_2>0$ such that~\cite{ocone1996}:
	\begin{align*}
	\|\tilde{\Sigma}_t-\Sigma_\infty\|_2 &\leq M_1e^{-2\lambda t}\\
	\Expect[|\tilde{m}_t-m_t|^2] &\leq M_2e^{-2\lambda t},	
	\end{align*}
	Explicit estimates of the constants $M_1,M_2$ appear in a recent paper~\cite{bishop2017stability}. 
	\subsection{Square root Riccati flow}
For the linear Gaussian filtering problem~\eqref{eq:dyn}-\eqref{eq:obs} define 
%
%
\begin{align*}
\sRicc (Q) &:= A - \frac{1}{2} Q H^\top H
\end{align*}
for $Q\in S^d_{++}$.
Also for $Q \in C([0,\infty),S^d_{++})$,  let $\Psi^{(Q)}_{t,s}$ be the state transition matrix for
\begin{equation}
\frac{\ud}{\ud t}\Psi^{(Q)}_{t,s} = \sRicc(Q_t)\Psi^{(Q)}_{t,s},\quad \Psi^{(Q)}_{s,s}=I\label{eq:sRiccati-flow}
\end{equation}	

The following Lemma is analogue of Lemma~\ref{lem:Riccati} for the square root Riccati flow~\eqref{eq:sRiccati-flow}. It is used to prove the stability of linear FPF~\eqref{eq:mean-field-FPF} in Prop.~\ref{prop:mean-field-stability}. 
The proof of Lemma~\ref{lem:sRiccati} appears in~Appendix~\ref{apdx:stability-KF}.
\begin{lemma}\label{lem:sRiccati}
Consider the linear flow~\eqref{eq:sRiccati-flow} with $Q_t=\Sigma_t$ where $\Sigma_t$ solves the DRE~\eqref{eq:DRE}. The state transition matrix $\Psi^{(\Sigma)}_{t,s}$ satisfies the bound
 		\begin{equation}\label{eq:Psi-bound}
 		\|\Psi^{(\Sigma)}_{t,s}\|_2 \leq \alpha e^{-\beta (t-s)},\quad \forall t\geq s>0 
 		\end{equation}
 		where
$\beta = \frac{\lambda_\text{min}(\Sigma_B)}{2\lambda_\text{max}(\Sigma_\infty)}$ and 
 $ \alpha=e^{\frac{\sqrt{\text{cond}(\Sigma_\infty)} M_1\|H^\top H\|}{2\beta}}\sqrt{\text{cond}(\Sigma_\infty)}$. 
\end{lemma}



%
 \medskip
 

\section{Analysis of the mean-field system}\label{sec:mean-field-process}
\subsection{Exactness}
 
Consider the mean-field sde~\eqref{eq:mean-field-FPF} for the FPF. 
Define the error process $\bar{\xi}_t:=\bar{X}_t-\bar{m}_t$. The evolution of the conditional mean $\bar{m}_t$, the conditional covariance  $\bar{\Sigma}_t$, and the  the error process $\bar{\xi}_t$ are given by the respective sdes:
\begin{subequations}
	\begin{align}
	\ud \bar{m}_t & = A \bar{m}_t \ud t + \bar{\Sigma}_t H^\top(\ud Z_t - H \bar{m}_t\ud t)\label{eq:mean-field-mean}\\\
	\frac{\ud}{\ud t} \bar{\Sigma}_t &= \Ricc(\bar{\Sigma}_t)\label{eq:DRE-mean-field}\\
	\ud \bar{\xi}_t &= \sRicc(\bar{\Sigma}_t)\bar{\xi}_t\ud t + \sigma_B \ud \bar{B}_t \label{eq:mean-field-error}
	\end{align}
\end{subequations}
%
Note that the sdes for the mean and covariance~\eqref{eq:mean-field-mean}-\eqref{eq:DRE-mean-field} are identical to the Kalman filter equations~\eqref{eq:KF-mean}-\eqref{eq:KF-variance}, and this property holds even for non-Gaussian initial distribution for $\bar{X}_0$.
\medskip

\begin{proposition}{(Exactness \cite[Thm. 1]{AmirACC2016})} \label{prop:exactness}
Consider the linear Gaussian filtering problem~\eqref{eq:dyn}-\eqref{eq:obs}, the Kalman-Bucy filter~\eqref{eq:KF-mean}-\eqref{eq:KF-variance}, and the McKean-Vlasov sde~\eqref{eq:mean-field-FPF} with initial distribution $\bar{X}_0\sim \bar{\pi}_0$. 
\begin{romannum}
	\item If $\bar{m}_0=m_0$ and $\bar{\Sigma}_0=\Sigma_0$ then
	\begin{equation*}
	\bar{m}_t=m_t,\quad \bar{\Sigma}_t = \Sigma_t,\quad \forall t>0
	\end{equation*}
	\item If the initial distribution $\bar{\pi}_0$ is Gaussian $\mathcal{N}(m_0,\Sigma_0)$ then 
	\begin{align*}
		\mathsf{P}(\bar{X}_t|\clZ_t) \text{ is Gaussian } \mathcal{N}(m_t,\Sigma_t)
	\end{align*}
\end{romannum}

\end{proposition}
\medskip

\begin{remark}
According to the Prop.~\ref{prop:exactness} the sde~\eqref{eq:mean-field-FPF} is exact for the case of Gaussian prior, i.e the conditional distribution of $\bar{X}_t$ is equal to the conditional distribution of hidden state $X_t$. The sde~\eqref{eq:mean-field-FPF} is not the only sde that satisfies the exactness property. In fact any sde of the following form is exact:
\begin{align*}
\ud \bar{X}_t = &A \bar{X}_t \ud t + \gamma_1 \sigma_B\ud \bar{B}_t + \frac{1-\gamma_1^2}{2}\bar{\Sigma}_t^{-1}(\bar{X}_t-\bar{m}_t)\ud t\\
&+ \bar{\k}_t(\ud Z_t - \frac{(1-\gamma_2^2)\bar{m}_t + (1+\gamma_2^2)\bar{X}_t}{2}\ud t+  \gamma_2\ud \bar{W}_t)  
\end{align*} 
where $\bar{W}_t$ is an independent copy of the observation noise, and $\gamma_1,\gamma_2$ are constants. The following three choices for $\gamma_1,\gamma_2$ lead to the three established forms of the linear FPF and EnKBF:
\begin{romannum}
	\item $\gamma_1=1$ and $\gamma_2=1$: EnKF with perturbed observation~\cite{jdw:BergemannReich2012}~\cite{delmoral2016stability};
	\item $\gamma_1=1$ and $\gamma_2=0$: Stochastic linear FPF~\cite[Eq. (26)]{yang2016} or square root EnKBF~\cite{jdw:BergemannReich2012} 
	\item $\gamma_1=0$ and $\gamma_2=0$: Deterministic linear FPF~\cite[Eq. (15)]{AmirACC2016}.
\end{romannum} 

\end{remark}
	
%
%
%
%

\subsection{Existence and uniqueness} 

To prove the existence of a unique solution $\bar{X}_t$ to the mean-field model~\eqref{eq:mean-field-FPF}, we use the decomposition $\bar{X}_t=\bar{m}_t+\bar{\xi}_t$ and we only consider the sde~\eqref{eq:mean-field-error} for the error process.
Given the solution $\bar{\xi}_t$, 
 the existence of the solution $\bar{m}_t$ of the sde~\eqref{eq:mean-field-mean} is straightforward 
because it is a linear sde. 

%
The proof for the existence relies on a fixed-point iteration and contraction argument. In order to state the result, the following definitions are necessary:
Let $T>0$ be the terminal time. For two random processes $X,Y$ on $C([0,T],\Re^d)$ define the metric 
\begin{equation*}
L_{2,t}(X,Y) := (\Expect[\sup_{s\in[0,t]} |X_s-Y_s|^2])^{\frac{1}{2}}
\end{equation*}
And let $\mathcal{X}$ be the space of random processes $X$ in $C([0,T],\Re^d)$ such that $L_{2,T} (X,0)<\infty$. 
For two positive definite symmetric valued functions $Q,U \in C([0,T],S^d_{++})$ define the metric \[l_{2,t}(Q,U) := \sup_{s \in [0,t]} \|Q_s-U_s\|_2 \] and let $\fSpace:=\{Q \in C([0,T],S^d_{++});~l_{2,t}(Q,0)<\infty\}$ . 

For a given $Q \in \mathcal{S}$ define the linear sde 
\begin{equation}
\ud \bar{\xi}^{(Q)}_t = \sRicc(Q_t) \bar{\xi}^{(Q)}_t\ud t + \sigma_B \ud \bar{B}_t, \label{eq:mean-field-error-Q}
\end{equation}
with $\bar{\xi}^{(Q)}_0 \sim \mathcal{N}(0,\Sigma_0)$.  Define the map $F:\mathcal{S}\to \mathcal{S}$ according to  $F(Q)_t:=\Expect[\bar{\xi}_t^{(Q)}(\bar{\xi}^{(Q)})^\top]$ for $t \in [0,T]$. Note that sde~\eqref{eq:mean-field-error} is of the form~\eqref{eq:mean-field-error-Q} with $Q_t=\bar{\Sigma}_t$ and $\bar{\Sigma_t}$ is a fixed point of the map $F$. 
The proof of the following Proposition appears in Appendix~\ref{apdx:existence}.
\medskip

\begin{proposition}[Existence of the mean-field process] \label{prop:mean-field-existence}
	Consider 
	the Mckean-Vlasov sde~\eqref{eq:mean-field-error} and the sde~\eqref{eq:mean-field-error-Q} for a fixed terminal time $T>0$.
	\begin{romannum}
		\item The sde~\eqref{eq:mean-field-error-Q} has a unique strong solution in $\mathcal{X}$. The map $F$ is well-defined and   satisfies the bound
	\begin{align*}
	l_{2,t}(F(Q)_t,0)&\leq
	e^{2\kappa t} \|\Sigma_0\|_2+ \frac{e^{2\kappa t}-1}{2\kappa}\|\Sigma_B\|_2
	\end{align*}	
	where $\kappa = \|A\|_2+\frac{l_{2,t}(Q,0)\|H^\top H\|_2}{2}$. 
		
		\item For  $Q,U\in \fSpace$ and $t<T$:
		\begin{subequations}
		\begin{align}
			l_{2,t}(F(Q),F(U)) &\leq C_2\int_{0}^t l_{2,s}(Q,U)\ud s\label{eq:F-Q-bound}\\
	L_{2,t}(\bar{\xi}^{(Q)},\bar{\xi}^{(U)}) &\leq C_1\int_{0}^t l_{2,s}(Q,U)\ud s\label{eq:xi-Q-bound}
	\end{align}			
		\end{subequations}
		where $C_1=e^{2\kappa t}\|H^\top H\|l_{2,t}(F(Q),0)$ and $C_2=\frac{\|H^\top H\|^2}{2}e^{\kappa t}l_{2,t}(F(Q),0)$.

\item 
The Mckean-Vlasov sde~\eqref{eq:mean-field-error} has a unique strong solution in~$\mathcal{X}$.
	\end{romannum} 
\end{proposition}
\medskip

\begin{remark}
The interaction terms in the stochastic linear FPF only effect the drift function, whereas the interaction in the EnKBF with perturbed observation and deterministic linear FPF also effect the process noise.
This difference makes the analysis of the stochastic linear FPF simpler.
\end{remark}

\subsection{Stability}
To prove the stability of the mean-field process $\bar{X}_t$,
consider the decomposition $\bar{X}_t=\bar{m}_t+\bar{\xi}_t$. The conditional mean $\bar{m}_t$ evolves according to~\eqref{eq:mean-field-mean} whose stability follows from Lemma~\ref{lem:Riccati}. The error process $\bar{\xi}_t$ evolves according to~\eqref{eq:mean-field-error} whose stability follows form Lemma~\ref{lem:sRiccati}. 
Therefore one may conclude that the mean-field process $\bar{X}_t=\bar{m}_t+\bar{\xi}_t$ is stable. The precise statement of the result is the following Proposition. The proof appears in Appendix~\ref{apdx:mean-field-stability}.  
\begin{proposition}(Stability of the mean-field process)\label{prop:mean-field-stability}
	Let $\bar{X}_t$ denote the solution to the McKean-Vlasov sde~\eqref{eq:mean-field-FPF} with the correct initial distribution $\mathcal{N}(m_0,\Sigma_0)$, and let $\tilde{X}_t$ denote the solution to the McKean-Vlasov sde~\eqref{eq:mean-field-FPF} with the initial distribution~$\tilde{\pi}_0$ that has finite second moment. Let $\bar{\pi}_t,\tilde{\pi}_t$ denote the conditional probability distribution of $\bar{X}_t$ and $\tilde{X}_t$ given $\clZ_t$ respectively. Then for all $t>0$, 
	\begin{equation*}
	W_2(\bar{\pi}_t,\tilde{\pi}_t) \leq  e^{-\beta t}\left( \alpha W_2(\pi_0,\tilde{\pi}_0) +M_2+\alpha M_1\|H^\top H\|_2\int_0^t \sqrt{\tr(\Sigma_s)}\ud s\right)
	\end{equation*} 
	where $M_1,M_2,\alpha$ are constants defined in Sec.~\ref{sec:stability-KF} and Lemma~\ref{lem:sRiccati}.
\end{proposition}
\medskip

\begin{remark}(Comparison with EnKBF)
	The Prop.~\ref{prop:mean-field-stability} is analogue to the result~\cite[Theorem 3.4]{delmoral2016stability} for the stability of the mean-field limit of the EnKBF with perturbed observation. The dynamics of the EnKBF with perturbed observation is governed by time-varying matrix $A-\Sigma_tH^\top H$ which has stronger stability properties compared to $A-\frac{1}{2}\Sigma_tH^\top H$ that governs the dynamics of stochastic linear FPF.     
\end{remark}
\section{Analysis of the finite-$N$ system } \label{sec:finite-N}
Consider the finite-$N$ system~\eqref{eq:finite-N-FPF}. 
Define the error processes $\xi^i_t:=X^i_t-\mN_t$ for $i=1,\ldots,N$. The evolutions for the empirical mean $\mN_t$, the empirical covariance $\SigN_t$, and the error $\xi^i_t$ are as follows
\begin{subequations}
	\begin{align}
	\ud \mN_t & = A \mN_t \ud t + \sigma_B\ud {B}^{(N)}_t + \k^{(N)}(\ud Z_t - H\mN_t\ud t) \label{eq:finite-N-mean}
	\\
	\ud \SigN_t &= \Ricc(\SigN_t)\ud t+ \ud M_t + \ud M_t^\top  \label{eq:finite-N-DRE}\\
	\ud \xi^i_t &= \sRicc(\SigN_t) \xi^i_t\ud t + \sigma_B \ud B^i_t -\sigma_B\ud {B}^{(N)}_t  \label{eq:finite-N-error}
	\end{align} 
\end{subequations}
where 
$\ud {B}^{(N)}_t := \frac{1}{N}\sum_{i=1}^N \ud B^i_t$ and $M^{(N)}_t$ is a matrix valued martingale given by
$
\ud M_t
:= \frac{1}{N-1} \sum_{i=1}^N
\sigma_B \ud B^i_t {\xi^i_t}^\top 
$.

The equations for the empirical mean~\eqref{eq:finite-N-mean} and the empirical covariance~\eqref{eq:finite-N-DRE} are similar to the Kalman filter equations~\eqref{eq:KF-mean}-\eqref{eq:KF-variance} except the additional stochastic terms $B_t^{(N)}$ and $M_t$ that scale as $O(\frac{1}{\sqrt{N}})$. 



We restrict the analysis to the scalar case ($d=1$). In addition to Assumption  A1 and A2, we make the following assumption:

\newP{Assumption A3} The state $X_t$ is asymptotically stable, i.e 
$\mu(A):=\min\{-\text{real}(\lambda);~\lambda~\text{is eigenvalue of}~A\} >0$. 

\subsection{Convergence of the empirical mean and covariance}
The main result regarding the convergence of the empirical mean and empirical covariance is the following Proposition. The proof appears in the Appendix~\ref{apdx:mean-var-conv}.
\begin{proposition}\label{prop:mean-var-conv}
	Consider the mean-field system~\eqref{eq:mean-field-FPF}, and the finite-$N$ system~\eqref{eq:finite-N-FPF} for the scalar case ($d=1$). 
	\begin{romannum}
		\item  For any $t>0$, and $N> 4p$:
			\begin{align}
			\Expect[|\Sigma_t^{(N)}-\Sigma_t|^{2p}]^{\frac{1}{p}} &\leq \frac{C_1}{N} e^{-2\beta t} + \frac{C_2}{N}  \label{eq:conv-var}
			\end{align}
		where $C_1=2\alpha^4 \Sigma_0^2 [(2p-1)!!]^{1/p}$, $C_2= 4(2p-1)\alpha^4\Sigma_\infty(\Sigma_0+\Sigma_\infty)$
		with $\alpha,\beta$  defined in Lemma~\ref{lem:sRiccati}, 
		\item  For any $t>0$ and as $N\to \infty $:
			\begin{align}
			\Expect[|\mN_t-m_t|^2] &\leq  \frac{\Sigma_0}{N}e^{-2\mu(A)t}+\frac{C_3}{N}\label{eq:conv-mean}
			\end{align}			
		where the constant
		$C_3= \frac{(C_1+C_2)H^2+\Sigma_B}{2\mu(A)}$.
	\end{romannum}
	\end{proposition}
	\medskip
	
\begin{remark}
	The result regarding the convergence of the empirical covariance~\eqref{eq:conv-var} follows without Assumption A3. Assumption A3 is required to prove the estimate~\eqref{eq:conv-mean}. 
\end{remark}


\subsection{Propagation of chaos analysis}
The next  objective to prove the convergence of the empirical distribution of the particles~$\{X^i_t\}_{i=1}^N$ to the distribution of the mean-field process $\bar{X}_t$, i.e,
\begin{equation*}
\frac{1}{N}\sum_{i=1}^N f(X^i_t) \to \Expect[f(\bar{X}_t)|\clZ_t]
\end{equation*}
for all bounded functions $f:\Re^d \to \Re$.

To show the convergence, introduce $N$ independent copies of the mean-field process $\bar{X}_t$
denoted by $\{\bar{X}_t^i;~i=1,\ldots,N\}$ such that 
\begin{equation} \label{eq:mean-field-copy}
\ud \bar{X}^i_t = A \bar{X}^i_t \ud t + \ud B^i_t + \bar{\k}_t (\ud Z_t - \frac{H\bar{X}^i_t+H\bar{m}_t}{2}\ud t),\quad \bar{X}^i_0=X^i_0
\end{equation}  
for $i=1,\ldots,N$. Note that $\bar{X}^i_t$ and $X^i_t$ are coupled through the same initial condition and the same process noise $\ud B^i_t$. 
Also define  the error process $\bar{\xi}^i_t:=\bar{X}^i_t-\bar{m}_t$. 
The error processes $\bar{\xi}^i_t$ evolve according to:
\begin{equation}\label{eq:mean-field-error-copy}
\ud \bar{\xi}^i_t = \sRicc(\bar{\Sigma_t})\bar{\xi}^i_t \ud t + \sigma_B \ud B^i_t
\end{equation}

\medskip
The result regarding the convergence of the empirical distribution is the following Proposition. The proof appears in Appendix~\ref{apdx:conv-POC}.
\begin{proposition}\label{prop:conv-POC}
	Consider the mean-field system~\eqref{eq:mean-field-FPF}, the finite-$N$ system~\eqref{eq:finite-N-FPF}, and the stochastic processes $\bar{X}^i_t$ defined in~\eqref{eq:mean-field-copy} for the scalar case ($d=1$). 
	\begin{romannum}
%
\item Particles: For any $t>0$ and as $N\to \infty $:
\begin{equation}\label{eq:conv-Xi}
\Expect[|X^i_t-\bar{X}^i_t|^2] \leq \frac{C_4}{N}
\end{equation}
for $i=1,\ldots,N$ where the constant $C_4= 2C_3 + 4\Sigma_0+\frac{\sqrt{3}H^4(\Sigma_0+\Sigma_\infty)(C_1+C_2)}{\mu(A)^2}+\frac{2\Sigma_B}{\mu(A)} $.
\item For any Lipschitz function $f$ 
\begin{equation}
\Expect\left[\left|\frac{1}{N}\sum_{i=1}^N f(X^i_t)-\Expect[f(\bar{X}_t)|\clZ_t]\right|^2\right] \leq \frac{(\text{const})}{N}\label{eq:conv-POC}
\end{equation}

	\end{romannum}
	
\end{proposition}
\medskip

\begin{remark}
Similar results for the vector case for the EnKBF with perturbed observation is shown in~\cite{delmoral2016stability}. 
The Assumptions in~\cite{delmoral2016stability} are (i) The matrix $A$ is stable; (ii) The matrix $H^\top H = \rho I$ (full rank observation matrix). 
\end{remark}

\medskip

\bibliographystyle{plain}
\bibliography{../bibfiles/fpfbib,../bibfiles/ref,../bibfiles/fpfbib2,../bibfiles/Optimization,../bibfiles/meanfield,../bibfiles/meanfield_v2} 

\section{Appendix}
\subsection{Proof of Lemma~\ref{lem:sRiccati}} \label{apdx:stability-KF}
	 	The proof is presented in two steps.
	
	\begin{enumerate}
		\item Consider the system  	  
		\begin{align*}
		\frac{\ud y_t}{\ud t } &= (\sRicc(\Sigma_\infty))^\top y_t
		\end{align*}	 
		with the Lyapunov function $V(y)=y^\top \Sigma_\infty y$. Observe
		\begin{align*}
		\frac{\ud }{\ud t} y_t^\top \Sigma_\infty y_t &= y_t^\top(A\Sigma_\infty + \Sigma_\infty A^\top -\Sigma_\infty H^\top H \Sigma_\infty)y_t\\
		&=-y_t^\top \Sigma_B y_t \leq -\frac{\lambda_\text{min}(\Sigma_B)}{\lambda_\text{max}(\Sigma_\infty)}  y_t^\top\Sigma_\infty y_t 
		\end{align*}
		where ARE~\eqref{eq:ARE} is used. Dividing both sides by $y_t\Sigma_\infty y_t$ and integrating with time yields:
		\begin{equation*}
		y_t\Sigma_\infty y_t \leq e^{-2\beta t} y_0\Sigma_\infty y_0 
		\end{equation*}
		where $\beta = \frac{\lambda_\text{min}(\Sigma_B)}{2\lambda_\text{max}(\Sigma_\infty)}$. 
		Using the inequality $\lambda_\text{min}(\Sigma_\infty)|y|^2\leq y^\top\Sigma_\infty y \leq \lambda_\text{max}(\Sigma_\infty)|y|^2$ yields the inequality
		\begin{equation*}
		|y_t|^2\leq e^{-2\beta t} \text{cond}(\Sigma_\infty)|y_0|^2
		\end{equation*}
		Using $|y_t|=|e^{t\sRicc(\Sigma_\infty)^\top}y_0|$ concludes 
		 \[\|e^{\sRicc(\Sigma_\infty)t}\|_2\leq \sqrt{\text{cond}(\Sigma_\infty)} e^{-\beta t}\]
		\item Consider the system 
		\begin{align*}
		\frac{\ud x_t}{\ud t } &= \sRicc(\Sigma_t)x_t\\
		&=\sRicc(\Sigma_\infty) x_t + \frac{1}{2}(\Sigma_\infty-\Sigma_t)H^\top Hx_t
		\end{align*}
		The solution satisfies the identity
		\begin{align*}
		x_t = &e^{\sRicc(\Sigma_\infty)t}x_0 \\&+ \frac{1}{2}\int_0^t e^{\sRicc(\Sigma_\infty)(t-s)}(\Sigma_\infty - \Sigma_s)H^\top H x_s \ud s
		\end{align*}
		and hence the inequality
		\begin{align*}
		|x_t|&\leq  \|e^{t\sRicc(\Sigma_\infty)}\|_2|x_0|\\&
		+\frac{\|H^\top H\|}{2}\int_0^t \|e^{(t-s)\sRicc(\Sigma_\infty)}\|_2 \|\Sigma_s-\Sigma_\infty\|_2 |x_s|\ud s 
		\end{align*} 
		Then use the inequality  from step~1,  $\|\Sigma_s-\Sigma_\infty\| \leq M_1 e^{-\beta s}$ from conclusion of Lemma~\ref{lem:Riccati} (because $\beta < \lambda_0$) 
to conclude 
		\begin{align*}
		|x_t|\leq  &e^{-\beta t}\sqrt{\text{cond}(\Sigma_\infty)}|x_0|
		\\&+\sqrt{\text{cond}(\Sigma_\infty)} \frac{M_1\|H^\top H\|}{2}\int_0^t e^{-\beta t}|x_s|\ud s 
		\end{align*} 	  
		Applying
  Gr\"onwall's inequality to $e^{\beta t} |x_t|$ yields:
		\begin{align*}
		|x_t|&\leq e^{-\beta t} e^{\sqrt{\text{cond}(\Sigma_\infty)} \frac{M_1\|H^\top H\|}{2}\int_0^te^{-\beta s} \ud s} \sqrt{\text{cond}(\Sigma_\infty)}|x_0|    \\&\leq  e^{-\beta t}e^{\frac{\sqrt{\text{cond}(\Sigma_\infty)} M_1\|H^\top H\|}{2\beta}}\sqrt{\text{cond}(\Sigma_\infty)}|x_0|
		\end{align*} 
		which proves the estimate~\eqref{eq:Psi-bound} for $s=0$. The result is extended to any time $s<t$ by following the same argument with the initial time $0$ replaced by $s$.
	\end{enumerate}



\subsection{Proof of the Proposition~\ref{prop:mean-field-existence}}\label{apdx:existence}
	\begin{romannum}
		\item A unique solution exists because the sde~\eqref{eq:mean-field-error-Q} is linear. The bound follows from writing down the solution explicitly,
		\[
		\bar{\xi}^{(Q)}_t = \Psi_t^{(Q)} \bar{\xi}_0 +  \int_0^t \Psi_{t,\tau}^{(Q)}  \sigma_B \ud \bar{B}_\tau
		\]  
		and the upper-bound $\|\Psi^{(Q)}_{t,\tau}\| \leq e^{\kappa(t-\tau)}$ with $\kappa = \|A\|+\frac{L_{2,t}(Q,0)\|H^\top H\|_2}{2}$. 
		\item For $Q,U \in \fSpace$ we have
		\begin{align*}
		\frac{\ud}{\ud t}(F(Q)_t - F(U)_t) =& \sRicc(Q_t)(F(Q)_t - F(U)_t)\\
		&+(F(Q)_t - F(U)_t)\sRicc(Q_t)^\top\\
		&+\frac{1}{2}(Q_t-U_t)H^\top H  F(U)_t  
		\\
		&+\frac{1}{2} F(U)_tH^\top H(Q_t-U_t) 
		\end{align*} 
		Therefore
		\begin{align*}
		F(Q)_t &- F(U)_t = \\&\frac{1}{2}\int_0^t \Psi^{(Q)}_{t,\tau}(Q_\tau - U_\tau)H^\top HF(U)_\tau  {\Psi_{t,\tau}^{(Q)}}^\top\ud \tau  \\
		+ &\frac{1}{2}\int_0^t \Psi^{(Q)}_{t,\tau}F(U)_\tau H^\top H(Q_\tau-U_\tau )  {\Psi_{t,\tau}^{(Q)}}^\top\ud \tau
		\end{align*} 	
		which satisfies the bound
		\begin{align*}
		l_{2,t}(F(Q),F(U)) 
		&\leq C\int_0^t l_{2,\tau}(Q,U)\ud \tau
		\end{align*}
		where the constant $C=e^{2\kappa t}\|H^\top H\|l_{2,t}(F(Q),0)$. To derive the bound~\eqref{eq:xi-Q-bound} note that	
		\begin{align*}
		\ud (\bar{\xi}^{(Q)}_t - \bar{\xi}^U_t)  &= \sRicc(Q_t)(\bar{\xi}^{Q}_t-\bar{\xi}^{(U)}_t)\ud t \\
		& + \frac{1}{2}(Q_t-U_t)H^\top H \bar{\xi}^{(U)}_t \ud t 
		\end{align*}  
		whose solution is
		\begin{equation*}
		\bar{\xi}_t^{(Q)} - \bar{\xi}_t^{(U)} = \frac{1}{2}\int_0^t (Q_\tau-U_\tau)H^\top H \bar{\xi}^{(U)} \ud \tau 
		\end{equation*}
		Therefore
		\begin{align*}
		L_{2,t}(\xi^{(Q)},\xi^{(U)}) 
&\leq 
		C \int_0^t l_{2,\tau }( Q,U) \ud \tau 
		\end{align*}
		where the constant $C=\frac{\|H^\top H\|^2}{2}e^{\kappa t}l_{2,t}(F(Q),0)$.
		\item Let $R=2l_{2,T}(\bar{\Sigma},0)$. Define the subset $\fSpace_R:=\{Q \in \fSpace;~l_{2,t_R}(Q,0)\leq R\}$. 
		According to the bound~\eqref{eq:F-Q-bound}, there exists $t_R>0$ small enough such that $f(\fSpace_R) \subset \fSpace_R$ and the map $F:\fSpace_R \to \fSpace_R$ is a contraction, i.e  
		\begin{equation*}
		l_2(F(Q),F(U)) < l_2(Q,U)
		\end{equation*}
		for all $Q,U\in \fSpace_R$. By definition, a random process $\bar{\xi}_t^{(Q)}$ is the solution of the sde~\eqref{eq:mean-field-error} iff $Q$ is the fixed point of the map $F$. Consider the fixed point iteration $Q^{(n+1)}=F(Q^{(n)})$ starting from $Q^{(0)}=I$.  
		According to the contraction mapping theorem the sequence converges to the fixed point. The fixed point of the map is equal to ${\Sigma}_t$ according to Proposition~\ref{prop:exactness}. The sequence $\bar{\xi}^{(Q^{(n)})}$ also converges due to the bound~\eqref{eq:xi-Q-bound}. As a result a unique strong solution exists on the interval $t\in [0,t_R]$. One may follow the same argument starting from $t_R$ instead of $0$ and extend the result further. This is possible because the fixed point always belongs to compact subset $\fSpace_R$.
	\end{romannum}

\subsection{Proof of the Proposition~\ref{prop:mean-field-stability}}\label{apdx:mean-field-stability}
\begin{proof}
Use the decomposition $\bar{X}_t = \bar{m}_t + \bar{\xi}_t$ and $\tilde{X}_t = \tilde{m}_t + \tilde{\xi}_t$. The conditional means $\bar{m}_t$, $\tilde{m}_t$ and the conditional covariances $\bar{\Sigma}_t$, $\tilde{\Sigma}_t$ evolve according to the Kalman filter equations~\eqref{eq:KF-mean}-\eqref{eq:KF-variance}. 
Therefore the difference $\Expect[|\bar{m}_t - \tilde{m}_t|^2]\leq M_2 e^{-2\lambda t}$ and $\|\bar{\Sigma}_t - \tilde{\Sigma}_t\|_2\leq M_1 e^{-2\lambda t}$ by Lemma~\ref{lem:Riccati}.

For the difference of the error processes $\bar{\xi}_t$ and $\tilde{\xi}_t$ we have
\begin{align*}
\ud (\tilde{\xi}_t-\bar{\xi}_t) = & \sRicc(\tilde{\Sigma}_t)(\tilde{\xi}_t-\bar{\xi}_t) \ud t+ \frac{1}{2}(\Sigma_t-\tilde{\Sigma}_t)H^\top H \bar{\xi}_t \ud t
\end{align*} 
The solution satisfies the identity
\begin{align*}
\tilde{\xi}_t-\bar{\xi}_t  = & \Psi^{(\tilde{\Sigma})}_t (\tilde{\xi}_0-\bar{\xi}_0)
+ \frac{1}{2}\int_0^t \Psi^{(\tilde{\Sigma})}_{t,s}(\Sigma_s-\tilde{\Sigma}_s)H^\top H \bar{\xi}_s \ud t
\end{align*}     
Use the result from Lemma~\ref{lem:sRiccati} to conclude the inequality
\begin{align*}
|\tilde{\xi}_t-\bar{\xi}_t&| \leq  \alpha e^{-\beta t} |\tilde{\xi}_0-\bar{\xi}_0|
\\&+ \alpha\frac{\|H^\top H\|_2}{2}\int_0^t e^{-\beta(t-s)}\|\Sigma_s-\tilde{\Sigma}_s\|_2 |\bar{\xi}_s| \ud t
\end{align*}   
Using the inequality $\|\Sigma_s-\tilde{\Sigma}_s\|_2\leq 2M_1e^{-\beta t}$ and taking  the mean-squred norm:
\begin{align*}
\Expect[\tilde{\xi}_t-\bar{\xi}_t|^2]^{1/2} \leq & \alpha e^{-\beta t}\Expect[|\tilde{\xi}_0-\bar{\xi}_0|^2]^{1/2}
\\&+\alpha M_1\|H^\top H\|_2e^{-\beta t}\int_0^t \Expect[|\bar{\xi}_s^2]^{1/2} \ud s
\end{align*}  
which together with the estimate $\Expect[|\bar{m}_t - \tilde{m}_t|^2]\leq M_2 e^{-2\lambda t}$ gives
\begin{align*}
\Expect[|\tilde{X}_t-\bar{X}_t|^2]^{1/2}&\leq \Expect[|\tilde{m}_t-\bar{m}_t|^2]^{1/2}+\Expect[|\tilde{\xi}_t-\bar{\xi}_t|^2]^{1/2}\\&\leq 
\alpha e^{-\beta t} \Expect[|\tilde{X}_0-\bar{X}_0|^2]^{1/2} + M_2e^{-2\lambda t}\\&+ 
\alpha M_1\|H^\top H\|_2e^{-\beta t}\int_0^t \sqrt{\tr(\Sigma_s)}\ud s
\end{align*}  
The result follows by definition of the $L^2$-Wasserstein distance, by taking the inf over all couplings between $\bar{X}_t$ and $\tilde{X}_t$.
\end{proof}
\subsection{Proof of the Proposition~\ref{prop:mean-var-conv}}\label{apdx:mean-var-conv}

%
\begin{romannum}  
	\item By subtracting~\eqref{eq:Riccati-flow} from~\eqref{eq:finite-N-DRE}, the evolution for the difference $\SigN_t - \Sigma_t$ is
	\begin{equation*}
	\ud (\SigN_t - \Sigma_t) = 2(A-\frac{\SigN_t+\Sigma_t}{2}H^2)(\SigN_t-\Sigma_t) + 2\ud M_t
	\end{equation*}
	Define $R_t=\Expect[(\SigN_t-\Sigma_t)^{2p}]$. By the application of It\^o's rule 
	\begin{align*}
	\frac{\ud R_t}{\ud t}  &=  \Expect[4p(A-\frac{\Sigma_t^{N} + \Sigma_t}{2}H^2)(\SigN_t-\Sigma_t)^{2p} ]+ 
	 \\&+
	p(2p-1){\frac{4\Sigma_B}{N-1}} \Expect[\Sigma_t^{(N)} (\SigN_t-\Sigma_t)^{2p-2}]
	\end{align*}
	Using the inequality $(A-\frac{\Sigma_t^{N} + \Sigma_t}{2}H^2) \leq \sRicc(\Sigma_t)$ and 	\[\SigN_t\leq (\SigN_t - \Sigma_t) + \Sigma_t \leq  \frac{1}{2\Sigma_\infty}(\SigN_t-\Sigma_t)^2+\Sigma_0+2\Sigma_\infty\]  and 
	$\Expect[|\SigN_t-\Sigma_t|^{2p-2}]\leq R_t^{\frac{p-1}{p}}$ yields:
	\begin{align*}
	\frac{\ud R_t}{\ud t}  &
	\leq  (4p\sRicc(\Sigma_t)+\frac{4p(2p-1)\Sigma_B}{2\Sigma_\infty (N-1)})R_t \\
	&+\frac{4p(2p-1)\Sigma_B}{N-1}(2\Sigma_\infty+\Sigma_0)R_t^{\frac{p-1}{p}}
	\end{align*}	
	Therefore,
	\begin{align*}
	\frac{\ud R_t^{1/p}}{\ud t} &\leq 4(\sRicc(\Sigma_t)+\frac{(2p-1)\Sigma_B}{2\Sigma_\infty (N-1)})R_t^{1/p} \\&+ \frac{4(2p-1)\Sigma_B(2\Sigma_\infty+\Sigma_0)}{N-1}
	\end{align*}
	Using the bound $e^{\int_0^t \sRicc(\Sigma_s)\ud s} \leq \alpha e^{-\beta t}$ from Lemma~\ref{lem:sRiccati} concludes 
	\begin{align*}
	R_t^{1/p} \leq &\alpha^4e^{-4\beta t +\frac{4(2p-1)\beta}{N-1}t}R_0^{1/p} \\+& \frac{\alpha^4}{4\beta - \frac{4(2p-1)\beta}{N-1}}{\frac{4(2p-1)\Sigma_B}{N-1}}(2\Sigma_\infty + \Sigma_0)
\end{align*}
Using $N>4p$, and $R_0^{1/p}\leq \frac{2\Sigma_0^2}{N}[(2p-1)!!]^{1/p}+O(\frac{1}{N^2})$ concludes the estimate~\eqref{eq:conv-var}.
	\item 
	The difference $\mN_t-m_t$ satisfies
	\begin{align*}
	\ud (\mN_t - &m_t) = (A-\SigN_tH^2)(\mN_t - m_t)\ud t +\ud S_t
	\end{align*}
	where $S_t$ is martingale given by
	\begin{equation*}
	\ud S_t = (\SigN_t-\Sigma_t)H(\ud Z_t-Hm_t\ud t) 
	+ \sigma_B \ud {B}^{(N)}_t
	\end{equation*}
	Therefore by application of It\"o rule
	\begin{align*}
	\frac{\ud}{\ud t} \Expect[(m_t-&\mN_t)^{2}]= \Expect[2(A-\SigN_tH^2)(\mN_t - m_t)^2] \\
	&+\Expect[(\SigN_t-\Sigma_t)^2H^2] + \frac{\Sigma_B}{N}
	\end{align*}
	Using Assumption A3 and the upper-bound~\eqref{eq:conv-var} for $p=1$
	\begin{align*}
	\frac{\ud }{\ud t} \Expect[(m_t-&\mN_t)^{2}]\leq -2\mu(A)\Expect[(\mN_t - m_t)^2] \\
	&+\frac{C_1e^{-2\beta t}+C_2}{N}H^2 + \frac{\Sigma_B}{N} 
	\end{align*}	
	 Application of the Gr\"onwal inequality concludes the estimate~\eqref{eq:conv-mean}.
\end{romannum}
\subsection{Proof of the Proposition~\ref{prop:conv-POC}}\label{apdx:conv-POC}
\begin{romannum}
	\item The difference $(\xi^i_t - \bar{\xi}^i_t)$ satisfies the sde:
	\begin{align*}
	\ud (\xi^i_t - \bar{\xi}^i_t) 
	&= (A-\frac{1}{2}\SigN_tH^2)(\xi^i_t-\bar{\xi}^i_t)\ud t\\
	&+(\SigN_t-\bar{\Sigma}_t)\frac{H^2}{2}\bar{\xi}^i_t\ud t - \sigma_B \ud {B}^{(N)}_t
	\end{align*} 
	Therefore by application of the It\"o rule
	\begin{align*}
	\frac{\ud}{\ud t} \Expect[&(\xi^i_t - \bar{\xi}^i_t)^2] = \Expect[2(A-\frac{1}{2}\SigN_tH^2)(\xi^i_t-\bar{\xi}^i_t)^2]\\
	&+\Expect[(\SigN_t-\bar{\Sigma_t})H^2\bar{\xi}^i_t(\xi^i_t-\bar{\xi}^i_t)]+\frac{\Sigma_B}{N}
	\end{align*} 
	Using the inequality $A-\frac{1}{2}\SigN_tH^2\leq -\mu(A)$, and
	\begin{align*}
	(\SigN_t-\bar{\Sigma_t})&H^2\bar{\xi}^i_t(\xi^i_t-\bar{\xi}^i_t) \leq \frac{\mu(A)}{2}(\xi^i_t-\bar{\xi}^i_t)^2 \\+& \frac{H^4}{2\mu(A)}(\SigN_t-\bar{\Sigma_t})^2(\bar{\xi}_t^i)^2
	\end{align*}
	yields
	\begin{align*}
	\frac{\ud}{\ud t} &\Expect[(\xi^i_t - \bar{\xi}^i_t)^2] 
\leq - \mu(A)\Expect[(\xi^i_t-\bar{\xi}^i_t)^2]\\
	&+\frac{H^4}{2\mu(A)}\Expect[(\SigN_t-\bar{\Sigma_t})^2{\xi^i_t}^2]+ \frac{\Sigma_B}{N}
	\end{align*}
	Then use the Cauchy-Schwarz inequality inequality 
	\[\Expect[(\SigN_t-\Sigma_t)^2\bar{\xi}_t^2]\leq \Expect[(\SigN_t-\Sigma_t)^4]^{1/2}\Expect[\bar{\xi}_t^4]^{1/2}\]  and the upper-bound~\eqref{eq:conv-var} for $p=2$ to conclude:
	 \begin{align*}
	\frac{\ud}{\ud t} &\Expect[(\xi^i_t - \bar{\xi}^i_t)^2] 
\leq 
	-\mu(A)\Expect[(\xi^i_t-\bar{\xi}^i_t)^2]\\	&+\frac{H^4}{2\mu(A)}\frac{(C_1+C_2)}{N}3\Sigma_t+ \frac{\Sigma_B}{N}
	 \end{align*}
	 Using the Gr\"onwall inequality, and $\Expect[{\xi}^i_0-\bar{\xi}^i_0|^2]=\frac{\Sigma_0}{N}$ yields
	\begin{align*}
	\Expect[|\xi^i_t-\bar{\xi}^i_t]|^2] &\leq e^{-\mu(A)t}\frac{\Sigma_0}{N}
	\\&+ \frac{3H^4(\Sigma_0+\Sigma_\infty)}{2N\mu(A)^2}(C_1+C_2)+ \frac{\Sigma_B}{N\mu(A)} 
	\end{align*} 
Combining this result with the estimate~\eqref{eq:conv-mean}  and the inequality 
	\begin{align*}
	\Expect[|X^i_t-\bar{X}^i_t|^2] \leq 2\Expect[|\mN_t-\bar{m}_t|^2]+2\Expect[|\xi^i_t-\bar{\xi}^i_t|^2]
	\end{align*}
	concludes the estimate~\eqref{eq:conv-Xi}.
	\item 
	Note that
	\begin{align*}
	\frac{1}{N}\sum_{i=1}^N f(X^i_t)-\Expect[f(\bar{X}_t)|\clZ_t] = & \frac{1}{N}\sum_{i=1}^N (f(X^i_t)-f(\bar{X}^i_t)) \\+&\frac{1}{N}\sum_{i=1}^N f(\bar{X}^i_t)-\Expect[f(\bar{X}_t)|\clZ_t]
	\end{align*} 
	Taking the mean-squared norm and using the triangle inequality yields
	\begin{align*}
	\Expect&\left[|\frac{1}{N}\sum_{i=1}^N f(X^i_t)-\Expect[f(\bar{X}_t)|\clZ_t]|^2\right]^{1/2} \leq  \\&+\sum_{i=1}^N \text{Lip}(f)\Expect[|\bar{X}^i_t-X^i_t|^2]^{1/2}+\frac{\text{var}(f)}{\sqrt{N}}
	\end{align*}
	where we used the function $f$ is Lipschitz, and  $\bar{X}^i_t$ are i.i.d. Using the result of part (i) concludes the proof. 
\end{romannum}	
\end{document}